\documentclass{birkjour}
\usepackage{graphicx}
\usepackage{grffile}
\usepackage{longtable}
\usepackage{wrapfig}
\usepackage{rotating}
\usepackage[normalem]{ulem}
\usepackage{amsmath}
\usepackage{textcomp}
\usepackage{amssymb}
\usepackage{capt-of}
\newtheorem{theorem}{Theorem}[section]

\newtheorem{corollary}[theorem]{Corollary}

\theoremstyle{definition}
\author{Rafael Villarroel-Flores}
\address{Universidad Autónoma del Estado de Hidalgo\\ Carretera Pachuca-Tulancingo km. 4.5\\ Pachuca 42184 Hgo.\\ MEXICO}
\email{rafaelv@uaeh.edu.mx}
\date{\today}
\title{Group Orders That Imply a Nontrivial \(p\)-Core}
\begin{document}

\maketitle
\begin{abstract}
Given a prime number \(p\) and a natural number \(m\) not divided by \(p\), we propose the problem of finding the smallest number \(r_{0}\) such that for \(r\geq r_{0}\), every group \(G\) of order \(p^{r}m\) has a non-trivial normal \(p\)-subgroup. We prove that we can explicitly calculate the number \(r_{0}\) in the case where every group of order \(p^{r}m\) is solvable for all \(r\), and we obtain the value of \(r_{0}\) for a case where \(m\) is a product of two primes.
\end{abstract}

\section{Introduction}
\label{introduction}
Throughout this note, \(p\) will be a fixed prime number. We use \(O_{p}(G)\) to denote the \(p\)-core of \(G\), that is, its largest normal \(p\)-subgroup.

We propose the following optimization problem: Given a number \(m\) not divisible by \(p\), find the smallest \(r_{0}\) such that every group having order \(n =p^{r}m\), with \(r\geq r_{0}\), has a nontrivial \(p\)-core \(O_{p}(G)\). Denote such number \(r_{0}\) by \(\Lambda(p,m)\). In Theorem \ref{existence}, we will prove that \(\Lambda(p,m)\) is well-defined for any prime \(p\) and number \(m\) (with \(p\nmid m\)). In Theorem \ref{solvable} we explicitly determine the value of \(\Lambda(p,m)\) in the case that all groups whose order have the form \(p^{r}m\) are solvable (for example, if \(m\) is prime or if both \(p\) and \(m\) are odd). Finally, in Section \ref{example}, we calculate \(\Lambda(2,15)\), a case that is not covered by the previous theorem.

We remark that the motivation for this research came from the search for examples of finite groups \(G\) such that the Brown complex \({\mathcal S}_p({G})\) of nontrivial \(p\)-subgroups of \(G\) (see for example \cite{1987-webb-subgroup-complexes} for the definition and properties) is connected but not contractible. It is known that \({\mathcal S}_p({G})\) is contractible when \(G\) has a nontrivial normal \(p\)-subgroup, and Quillen conjectured in \cite{MR493916} that the converse is also true.

\section{Theorems}
\label{sec:org7220b64}

\begin{theorem}
\label{existence}
For any prime number \(p\) and natural number \(m\) such that \(p\nmid m\), there is a number \(\Lambda(p,m)\) such that if \(r\geq \Lambda(p,m)\), any group of order \(p^{r}m\) has a non-trivial \(p\)-core \(O_{p}(G)\).
\end{theorem}

\begin{proof}
Let \(G\) be a group of order \(p^{r}m\) with \(O_{p}(G)=1\). Let \(P\) be a Sylow \(p\)-subgroup of \(G\). Since the kernel of the action of \(G\) on the set of cosets of \(P\) is precisely \(O_{p}(G)\), we obtain that \(G\) embeds in \(S_{m}\), and so \(p^{r}\) divides \((m-1)!\). Hence, if \(p^{r_{0}}\) is the largest power of \(p\) dividing \((m-1)!\), we obtain that \(\Lambda(p,m)\leq r_{0}+1\).
\end{proof}

For \(t,q\) natural numbers, let \(\gamma(t,q)\) be the product
  \begin{equation}
\label{eq-gammatq}
\gamma(t,q)=(q^{t}-1)(q^{t-1}-1)\cdots (q^{2}-1)(q-1),
  \end{equation}
and if \(m=q_{1}^{t_1}q_{2}^{t_2}\cdots q_{k}^{t_k}\) is a prime factorization of \(m\), with the \(q_i\) pairwise distinct and \(t_i>0\) for each \(i\), we let \(\Gamma(m)=\gamma(t_{1},q_{1})\cdots\gamma(t_{k},q_{k})\). We prove that if \(p^{s_{0}}\) is the largest power of \(p\) dividing \(\Gamma(m)\), then \(\Lambda(p,m)\geq s_{0}+1\).

\begin{theorem}
  \label{lowerbound}
Let \(n=p^{s}m\) where \(p\nmid m\) and \(s>0\). If \(p^{s}\mid\Gamma(m)\), then there is a group of order \(n\) with \(O_{p}(G)=1\).
\end{theorem}

\begin{proof}
Let \(K\) be the group \(C_{q_1}^{t_1}\times\cdots\times
  C_{q_k}^{t_k}\), that is, a product of elementary abelian groups, where \(m=q_{1}^{t_1}\cdots q_{k}^{t_k}\) and \(q_{1},\ldots,q_{k}\) are distinct primes and \(C_{q}\) denotes the cyclic group of order \(q\). Then \(\Gamma(m)\) divides the order of \(\mathop{\mathrm{Aut}}(K)\), and hence  so does \(p^{s}\). Let \(H\) be a subgroup of \(\mathrm{Aut}(K)\) of  order \(p^{s}\). For every \(S\in H\) and \(k\in K\) define the map  \(T_{S,k}\colon K\to K\) by \(T_{S,k}(x)=Sx+k\). Then  \(G=\left\{\,T_{S,k}\mid S\in H, k\in K\,\right\}\) is also a subgroup of \(\mathrm{Aut}(K)\). If we identify \(H\) with the subgroup of maps of the form  \(T_{S,0}\) and \(K\) with the subgroup of maps of the form \(T_{1_K,k}\),  then \(G\) is just the semidirect product of \(K\) by \(H\). Hence \(|G|=n\). We  have that \(G\) acts transitively on \(K\) in a natural fashion, and the  stabilizer of \(0\in K\) is \(H\), a \(p\)-Sylow subgroup of \(G\). Hence the  stabilizers of points in \(K\) are precisely the Sylow subgroups of \(G\), so their intersection \(O_{p}(G)\) contains only the identity \(K\to K\), as we wanted to prove. 
\end{proof}

The next theorem will show that the lower bound given by Theorem \ref{lowerbound} is tight in some cases.

\begin{theorem}
\label{solvable}
Let \(n=p^{s}m\), where \(p\nmid m\). If \(G\) is a group of order \(n\) and \(p^{s}\) does not divide \(\Gamma(m)\) then either:
\begin{enumerate}
\item \(O_{p}(G)\ne 1\), or
\item \(G\) is not solvable.
\end{enumerate}
\end{theorem}

\begin{proof}
Let \(G\) be solvable with order \(n=p^{s}m\) and \(O_{p}(G)=1\). Let \(F(G)\) be the Fitting subgroup of \(G\). Consider the map \(c\colon G\to\mathop{\mathrm{Aut}}(F(G))\), sending \(g\) to \(c_{g}\colon F(G)\to F(G)\) given by conjugation by \(g\). The restriction of \(c\) to \(P\), a \(p\)-Sylow subgroup of \(G\), has kernel \(P\cap C_{G}(F(G))\). Since \(C_{G}(F(G))\le F(G)\) (Theorem 7.67 from \cite{MR1298629-djvu}), and \(F(G)\) does not contain elements of order \(p\) by our assumption on \(O_{p}(G)\), we have \(P\cap C_{G}(F(G))=1\) and so \(P\) acts faithfully on \(F(G)\). If \(m=q_{1}^{t_1}\cdots q_{k}^{t_k}\) is the prime factorization of \(m\), we have that \(F(G)\) is the direct product of the \(O_{q_{i}}(G)\) for \(i=1,\ldots,k\). Hence \(P\le\mathop{\mathrm{Aut}}(F(G))\cong \mathop{\mathrm{Aut}}(O_{q_1}(G))\times\cdots\times \mathop{\mathrm{Aut}}(O_{q_k}(G))\). Let \(g\in P\) such that the action induced by \(c_{g}\) on \(\prod_{i} O_{q_i}(G)/\Phi(O_{q_i}(G))\), is the identity. Since \(c_{g}\) acts on each factor \(O_{q_i}(G)/\Phi(O_{q_i}(G))\) as the identity, then by Theorem 5.1.4 from \cite{MR569209-big}, we have that it acts as the identity on each \(O_{q_i}(G)\). By the faithful action of \(P\) on \(F(G)\), we have that \(g=1\). This implies that \(P\) acts faithfully on \(\prod_{i} O_{q_i}(G)/\Phi(O_{q_i}(G))\). But then \(|P|\) divides the order of the automorphism group of \(\prod_{i} O_{q_i}(G)/\Phi(O_{q_i}(G))\), which is a product of elementary abelian groups of respective orders \(q_{i}^{s_i}\) with \(s_{i}\le t_{i}\) for all \(i\). Hence \(p^{s}=|P|\) divides \(\Gamma(m)\).
\end{proof}

\begin{corollary}
Let \(p^{s}\) be the largest power of \(p\) that divides \(\Gamma(m)\). If \(m\) is prime, or if both \(p,m\) are odd, then \(\Lambda(p,m)=s+1\).
\end{corollary}

\begin{proof}
By Burnside's \(p,q\)-theorem, and the Odd Order Theorem, we have that all groups that have order of the form \(p^{r}m\) for some \(r\) are solvable. Therefore, for all \(r>s\), by Theorem \ref{solvable} we have that all groups of order \(p^{r}m\) have non-trivial \(p\)-core.
\end{proof}

At this moment, we can prove that in some cases, the group constructed in \ref{lowerbound} is unique.

\begin{theorem}
  \label{unique}
Let \(n=p^{s}m\) where \(p\nmid m\) and \(s>0\). If \(p^{s}\mid\Gamma(m)\), but \(p^{s}\nmid\Gamma(m')\) for all proper divisors \(m'\) of \(m\), then up to isomorphism, the group constructed in the proof of Theorem \ref{lowerbound} is the only solvable group of order \(n\) with \(O_{p}(G)=1\).
\end{theorem}

\begin{proof}
With the notation of the argument of the proof of \ref{solvable}, if \(G\) is a solvable group of order \(n\) with \(O_{p}(G)=1\), we must have that \(|O_{q_{i}}(G)|=q_{i}^{t_{i}}\) and \(\Phi(O_{q_i}(G))=1\) for all \(i\) in order to satisfy the divisibility conditions. Hence \(O_{q_{i}}(G)\) is elementary abelian and a \(q_{i}\)-Sylow subgroup for all \(i\), and so \(G\) is the semidirect product of a \(p\)-Sylow subgroup \(P\) of \(F(G)=C_{q_1}^{t_1}\times\cdots\times C_{q_k}^{t_k}\) with \(F(G)\), where the action of \(P\) on \(F(G)\) by conjugation is faithful. Hence \(G\) is isomorphic to the group constructed in the proof of Theorem \ref{solvable}. 
\end{proof}

One case in that we may apply Theorem \ref{unique} is when \(n=864\). There are \(4725\) groups of order \(864=2^{5}3^{3}\), but only one of them has the property of having a trivial \(2\)-core.

\section{An example}
\label{sec:org022a6d6}
\label{example}

An example that cannot be tackled with the previous results is the case \(p=2\), \(m=3\cdot 5=15\). In this case, \(\Gamma(15)=(3-1)(5-1)=2^{3}\). Not all groups with order of the form \(2^{r}\cdot 3\cdot 5\) are solvable, however, we will prove that \(\Lambda(2,15)\) is actually \(4\). (The group \(S_{5}\) attests that \(\Lambda(2,15)>3\).)

\begin{theorem}
Every group \(G\) of order \(2^{r}\cdot 3\cdot 5\)  for \(r\geq 4\) is such that \(O_{2}(G)\ne 1\).
\end{theorem}

\begin{proof}
Let \(G\) be a group of order \(2^{r}\cdot 3\cdot 5\) for \(r\geq4\). Suppose that \(O_{2}(G)=1\). From Theorem \ref{solvable}, we obtain that \(G\) is not solvable. We will prove then that \(O_{3}(G)= 1\). Suppose otherwise, and let \(T=O_{3}(G)\). Then \(|G/T|=2^{r}\cdot 5\), and so \(G/T\) is solvable. Since \(2^{r}\nmid \Gamma(5)\), from Theorem \ref{solvable}, we have that \(O_{2}(G/T)\ne 1\). Let \(L\lhd G\) such that \(O_{2}(G/T)=L/T\). Suppose \(|L/T|=2^{j}\). Since \(O_{2}(G/L)=1\), \(|G/L|=2^{r-j}\cdot 5\) and \(G/L\) is solvable, we have that \(2^{r-j}\) divides \(\Gamma(5)=2^{2}\), that is, \(r-j\leq 2\). Now, \(L\) is also solvable and \(\Gamma(3)=3-1=2\), hence if we had \(j\geq 2\) we would have \(O_{2}(L)\ne 1\), and \(G\) would have a non-trivial subnormal \(2\)-subgroup, which contradicts our assumption that \(O_{2}(G)=1\). Hence \(j=1\). But then \(r-1\leq 2\), which contradicts that \(r\geq 4\). Hence \(O_{3}(G)=1\). By a similar argument, we get that \(O_{5}(G)=1\).

From \cite{1968-brauer-on-simple-groups-of-order-5dot-3-spadot-2-spb} we obtain that \(G\) is not simple. Hence \(G\) has a proper minimal normal subgroup \(M\). From the previous paragraph, we obtain that \(M\) is not abelian, since in that case we would have that \(M\leq F(G)\). The only possibility is that \(M=A_{5}\).  We have then a morphism \(c\colon G\to \mathrm{Aut}(A_{5})\)  sending \(g\) to \(c_{g}\), the conjugation by \(g\). Since \(\mathrm{Aut}(A_{5})=S_{5}\), and \(|c(G)|=|\mathrm{Inn}(G)|\geq |\mathrm{Inn}(A_{5})|=60\), in any case the kernel of \(c\) is a nontrivial normal \(2\)-subgroup.
\end{proof}

\bibliographystyle{plain}
\bibliography{../../../texmf/bibtex/bib/misc/rvf}
\end{document}